\documentclass[a4paper,11pt,reqno]{article}
\usepackage[T1]{fontenc}
\usepackage{courier}
\usepackage[scaled=0.92]{helvet}
\usepackage{amscd,amsmath,amssymb,amsthm,amsfonts,epsfig,graphics,accents}

\usepackage[percent]{overpic} 				
\usepackage[english]{babel}
\usepackage{algpseudocode,algorithm,float}
\usepackage{amsthm,amsfonts,amssymb,amsmath,subfigure,array,booktabs,longtable}
\usepackage{hyperref}
\hypersetup{
colorlinks=true,
citecolor=red,
linkbordercolor={1 1 1},
citebordercolor={1 1 1},
}
\usepackage{cleveref}
\usepackage{mathrsfs}
\usepackage{color}
\usepackage{bm}

\usepackage{xcolor} 

\newcommand{\mx}{\mbox}

\newcommand{\pl}{\partial}

\newcommand{\beq}[1]{\begin{equation} \label{#1}}
\newcommand{\eeq}{\end{equation}}
\newcommand{\beqar}{\[ \begin{array}{rcl}}
\newcommand{\eeqar}{\end{array} \]}

\newcommand{\ue}[1]{\boldsymbol{#1}}
\newcommand{\uue}[1]{{\boldsymbol{#1}}}

\newtheorem{rem}{Remark}[section]

\usepackage{xcolor}

\usepackage{xcolor}

\usepackage[utf8]{inputenc}

\hypersetup{pdfauthor={Fortunati, A. et al.},%
            pdftitle={A class of quasi-equatorial controlled switching
solutions to the iLLG equation},%
            pdfsubject={},%
            pdfkeywords={}%
}

\makeatletter
\renewcommand*{\@fnsymbol}[1]{\ensuremath{\ifcase#1\or *\or \mathsection \or \dagger \or \else \fi}}
\makeatother

\title{Quasi-equatorial controlled switching in inertial magnetization dynamics}
\date{}
\author{%
Alessandro Fortunati
$\,$ and Massimiliano d'Aquino\thanks{E-mail: mdaquino@unina.it}    
\bigskip \\
\textsuperscript{*}Department of Electrical Engineering and Information Technologies, \\University of Naples Federico II, I-80125 Naples, Italy.
}

\usepackage{amsmath}
\usepackage{graphicx}
\begin{document}

\maketitle

\begin{abstract}
We design exponentially decaying transverse magnetic fields
for reliable and ultra-fast (picosecond) magnetization reversal governed by the inertial
Landau--Lifshitz--Gilbert equation. The control construction
uses an explicitly resolvable auxiliary model to select the
field amplitude and characterize the influence of its decay
rate, with the aim of obtaining reversal followed by relaxation
without requiring a precisely timed interruption of the excitation. The auxiliary model is obtained through an equatorial
approximation, and bounds on its latitude identify a regime
consistent with quasi-equatorial motion.
For the full nonlinear
macrospin dynamics, an energy argument accounting for both
inertial kinetic energy and the remaining field input provides
a sufficient condition excluding backhopping, even while the
exponential tail remains applied. Approximate confinement times
and optimal decay rates are derived from the competition between
field decay and inertial relaxation. Numerical simulations
confirm the usefulness of these design estimates and identify
their limitations near a capture boundary, where nearby pulses
can lead to opposite final states. Spatially resolved inertial
micromagnetic simulations further support reversal followed by
relaxation over the simulated interval and quantify the transient
departures from spatially coherent motion. The resulting approach
combines explicit pulse design with a full-system energy
certificate, while distinguishing rapid initial rotation from
safe confinement and subsequent relaxation.

 
\end{abstract}

\section{Introduction}

Controlling magnetization reversal on short time scales is a central problem
in applied magnetism, motivated by the development of magnetic memories and
spintronic devices~\cite{DPG20}. Among the mechanisms available for this
purpose, precessional switching exploits the large-angle motion induced by
a magnetic field transverse to the initial magnetization, allowing reversal
without waiting for damping to drive the entire transition~\cite{KR02,BMS2009}. Its effectiveness depends, however, on the interplay
between the excitation, magnetic anisotropy and dissipation: reaching the
opposite hemisphere does not by itself ensure subsequent relaxation to the
desired equilibrium.

A characteristic difficulty of this approach is the sensitivity to the
duration of the applied field. When a constant transverse field is used,
its termination must place the magnetization in the basin of attraction of
the reversed equilibrium. Geometrical and analytical descriptions of this
process have clarified the role of the pulse duration and the subsequent
relaxation~\cite{BMSd03,dSS04,BMS2009}. In weakly dissipative regimes, the
dependence of the final state on the state at field removal can be
particularly intricate~\cite{SdBM09}. This makes the design of switching
protocols that do not require a precisely timed interruption of the
excitation an important control problem.

At sufficiently short time scales, magnetic inertia introduces an
additional feature into this picture. 

The inertial extension of the Landau--Lifshitz--Gilbert
equation (iLLG) incorporates a second-order time derivative
of the magnetization and predicts nutational motion~\cite{CRW11}. This inertial spin nutation occurs at terahertz frequencies in contrast to nutation-type magnetization motion that can also be produced in classical precessional dynamics at gigahertz frequencies~\cite{Serpico2004}.
Microscopic theoretical descriptions support the inclusion
of magnetic inertia on quantum mechanical grounds~\cite{Mondal2017,webb2026microscopic}, while experimental evidence for
nutational dynamics has been reported in ferromagnetic films~\cite{NAK21,Unikandanunni2022,de2025magnetic}.
The corresponding dynamics depends on
both the magnetization and its velocity, and its energy includes an
inertial kinetic contribution~\cite{NPS22,DAQUINO2024112874}. Consequently, the
magnetization direction alone is insufficient to characterize the state
from which relaxation begins. Numerical investigations have shown that
inertial effects can modify both the switching response and the associated
energy transfer~\cite{NPS22}. Nutational dynamics has also been exploited to induce switching
under high-frequency magnetic fields~\cite{Winter2022}, highlighting
magnetic inertia as an additional resource for controlling
magnetization reversal. 

An analytical treatment of field-driven switching in this regime was
developed in Ref.~\cite{FdS26}. Using a
multiple-time-scale perturbative approach, that work obtained approximate
solutions for the driven motion and explicit estimates of the switching
time and its admissible tolerance interval, with validation by macrospin
and micromagnetic simulations. The excitation considered there is a
constant field applied over a finite time interval. The resulting
estimates therefore quantify when the field may be removed while still
obtaining reversal. The present work addresses a complementary question:
can the temporal profile of the field itself be chosen so that successful
switching does not rely on selecting a separate switch-off instant?

To this end, we consider a macrospin initially at rest along its easy axis
and subject to a field directed along the third principal axis, with an
exponentially decreasing amplitude. We seek \emph{quasi-equatorial}
motions, in which the magnetization rotates towards the opposite easy-axis
orientation while its latitude remains small. This geometric objective
is related to the quasi-latitudinal control strategy introduced for
spin-torque-driven magnetization dynamics in Ref.~\cite{FdS24}; here the
control is supplied by an applied magnetic field and the dynamics includes
inertial effects. Rather than deriving a perturbative solution of the full
nonlinear equations, we construct the field using an explicitly resolvable
auxiliary system obtained by freezing the geometric coefficients at the
equator and neglecting the anisotropy torques.

The auxiliary dynamics suggests a family of controls of the form
\[
 h_3(t)=A\exp\!\left(-\frac{\gamma\alpha}{\eta}t\right),
 \qquad t\geq0,
\]
where $\alpha$ and $\eta$ denote the damping and inertial parameters,
respectively, and $\gamma>0$ adjusts the decay rate. The amplitude $A$ is
chosen so that the auxiliary azimuth approaches the reversed orientation.
The resulting formulas identify $\eta/\alpha$ as a natural decay time and
predict a limiting latitude of order $\alpha$. They also reveal a trade-off:
increasing $\gamma$ shortens the applied pulse but can increase transient
departures from the equator. Comparison with the full nonlinear iLLG
equation is therefore essential for assessing the effectiveness of the
constructed field.

A central feature of the proposed strategy is that the control
is designed explicitly around magnetic inertia. Indeed, its
amplitude and decay time are
\[
 A=\frac{\pi\gamma\alpha(1+\alpha^2)}{\eta},
 \qquad
 \tau_h=\frac{\eta}{\gamma\alpha},
\]
so the inertial parameter $\eta$ directly determines the
strength and temporal profile of the excitation. The field
decay is matched to the natural damping time $\eta/\alpha$
of the inertial velocity dynamics, while $\gamma$ adjusts
their relative time scales. Magnetic inertia is therefore
an essential ingredient of the control construction, rather
than a correction to a pulse designed within the classical
Landau--Lifshitz--Gilbert framework.

The main novelty is an inertia-based switching strategy in
which the field amplitude and decay rate are selected jointly
from the iLLG dynamics. The excitation decays according to a
prescribed law, avoiding the need to interrupt it at a
particular phase of the magnetization motion. For the parameter
regimes investigated numerically, this produces reversal
followed by relaxation towards the opposite easy-axis
equilibrium. The proposed strategy thus combines rapid
precessional motion with a timing tolerance resembling that
of the much slower 'damping' switching~\cite{BMS2009}. When field removal
is desired, a sufficient lower time threshold can be identified:
once energy confinement is established, the field may be
removed at any later time without compromising reversal.

The control-design problem is
thus expressed in terms of the pulse amplitude and decay rate, with
explicit auxiliary solutions providing guidance for their selection.
The principal formulas are collected in
Table~\ref{tab:control-design-formulas} as a practical guide
to pulse selection and assessment of safe confinement.

The paper is organized as follows. Section~\ref{sec:setup} introduces the spherical
coordinates and the auxiliary equatorial model. Section~\ref{sec:control} constructs
the exponentially decaying controls and derives their explicit
auxiliary solutions. Section~\ref{sec:relaxation} develops the energy confinement
criterion and analytical estimates for admissible decay rates,
safe switch-off times and optimal control selection. Section~\ref{sec:validation}
validates the approach through full nonlinear macrospin and inertial
micromagnetic simulations. Section~\ref{sec:conclusions} reports the conclusions.

\section{Set up and equatorial approximation}\label{sec:setup}
The ultra-fast magnetization dynamics of a uniformly magnetized (macrospin) particle is described by the inertial Landau--Lifshitz--Gilbert equation:
\beq{eq:illg}
\dot{\ue{m}}=-\ue{m} \times \bigl(\ue{h}_{\mathrm{eff}}-\alpha\,\dot{\ue{m}}-\eta\,\ddot{\ue{m}}\bigr) \mx{,}
\eeq
{where $\bm m$ is magnetization unit-vector (normalized by the saturation magnetization $M_s$), time is measured in units of $(\gamma_e M_s)^{-1}$ ($\gamma_e\approx\mu_0 |e|/m_e\approx 2.211\times 10^5\,\, m A^{-1} s^{-1}$ is the absolute value of the gyromagnetic ratio, being $e,m_e$ the electron charge and mass, respectively), $\alpha,\eta>0$ are the (dimensionless) Gilbert damping parameter and the strength of magnetic inertia, respectively (typically in the order of $10^{-2}$, see Ref.~\cite{daquino2023micromagnetic}). For a broader account of the modelling and numerical simulation
of inertial magnetization dynamics, see
Ref.~\cite{daquino2026micromagneticmodellingsimulationinertial}.
The effective magnetic field (also normalized by $M_s$) is}
\beq{eq:heff}
\ue{h}_{\mathrm{eff}}(t)=\ue{h}(t)-\bigl(D_1 m_1,\;D_2 m_2,\;D_3 m_3\bigr) \mx{,}
\eeq
where $\ue{h}(t)=(h_1(t),h_2(t),h_3(t))$ and
\beq{eq:dj}
-\infty < D_1 < D_2 < D_3 < +\infty 
\eeq
{are effective demagnetizing factors taking into account the (magnetostatic) shape and the magneto-crystalline anisotropy of the magnetic particle.
Inertial magnetization dynamics must fulfill the generalized free energy balance equation~\cite{DAQUINO2024112874}:
\begin{equation}\label{eq:iLLG energy balance}
    \dot{E}=-\alpha|\dot{\bm m}|^2-\dot{\bm h}\cdot\bm m \,,
\end{equation}
where the free energy $E$ (measured in units of $\mu_0M_s^2 V$, with $V$ being the particle volume) is given by:
\begin{equation}\label{eq:free energy}
    E(\bm m,\dot{\bm m},\bm h)=\frac{1}{2}\eta|\dot{\bm m}|^2+\frac{1}{2}(D_1m_1^2+D_2m_2^2+D_3m_3^2)-\bm h\cdot \bm m \,\,. 
\end{equation}}
We shall also denote $D_{i,j}:=D_i-D_j$ for all $i > j$.\\
We consider magnetization trajectories starting from the
easy-axis equilibrium at rest,
\[
 (\mathbf m(0),\dot{\mathbf m}(0))=(\mathbf e_1,\mathbf0).
\]
Then, we define the motion a successful magnetization switching driven by the applied field $\bm h(t)$ if
\[
 \lim_{t\to\infty}
 (\mathbf m(t),\dot{\mathbf m}(t))
 =(-\mathbf e_1,\mathbf0).
\]
Entry into the hemisphere $m_1<0$ is therefore only an
intermediate condition: successful reversal also requires
subsequent relaxation without a return to the original well.
A crossing from $m_1<0$ back to $m_1>0$ is hereafter termed backhopping.\\
Let us now consider the system of spherical coordinates 
\beq{eq:spherical}
m(t)=
\begin{pmatrix}
r(t)\cos\theta(t)\cos\phi(t)\\
r(t)\cos\theta(t)\sin\phi(t)\\
r(t)\sin\theta(t)
\end{pmatrix},
\eeq
where $r(t)=1$ and the equator is represented by $\theta=0$. In this reference, the iLLG equation (\ref{eq:illg}) reads as
\beq{eq:szero}
\begin{aligned}
\eta\Bigl(\cos\theta\,\ddot\phi-2\sin\theta\,\dot\theta\,\dot\phi\Bigr)
+\alpha\cos\theta\,\dot\phi
-\dot\theta
+\Bigl(D_{2,1} \cos\theta\cos\phi+h_1\Bigr)\sin\phi
-h_2\cos\phi
&=0.\\
\eta\Bigl(-\cos^{2}\theta\sin\theta\,\dot\phi^{2}-\cos\theta\,\ddot\theta\Bigr)
-\alpha\cos\theta\,\dot\theta
-\cos^{2}\theta\,\dot\phi
+D_{2,1} \cos^{2}\theta\sin\theta\,\sin^{2}\phi & \\
-h_2\cos\theta\sin\theta\,\sin\phi
-h_1\cos\theta\sin\theta\,\cos\phi
-D_{3,1} \cos^{2}\theta\sin\theta
+h_3\cos^{2}\theta
&=0.
\end{aligned} \tag{S}
\eeq
The latter are easily obtained by left multiplying each term of (\ref{eq:illg}) by $J^{\top}$, where $J:=\pl_{r,\theta,\phi}m$ is the Jacobian matrix of (\ref{eq:spherical}).\\
From now on, we shall consider the applied field of the form
\beq{eq:field}
\ue{h}(t)=(0,0,h_3(t)) \mx{.}
\eeq
As anticipated, the paper is concerned with the so-called \emph{equatorial approximation} i.e. the study of the system dynamics "near" $\theta=0$. Our aim is clearly to construct a control $h_3(t)$ such that the variations of $\theta(t)$ are bounded by a constant which is as small as possible.\\
Finally, recalling that typically for thin-film shaped magnetic particles $D_{2,1} =O(10^{-1})$, i.e. significantly smaller than $D_{3,1},D_{3,2}=O(1)$, we can assume it to be ''negligible'', this is an approximation which, according to the numerical experiments reported later on, is shown to be particularly adequate for this model.\\ 
{To construct an explicitly resolvable auxiliary system, we freeze
the geometric coefficients at the equator, replacing
$\sin\theta$ by $0$ and $\cos\theta$ by $1$, while retaining
$\theta(t)$ and its derivatives as dynamical variables.
We additionally omit the in-plane anisotropy torque by setting
$D_{2,1}=0$ in this auxiliary model.}

{The replacements defining the auxiliary model are therefore
\begin{equation}
 \sin\theta\mapsto0,\qquad
 \cos\theta\mapsto1,\qquad
 D_{2,1}\mapsto0,
 \label{eq:assumptions}
\end{equation}}
By using (\ref{eq:field}) and (\ref{eq:assumptions}), system \eqref{eq:szero} simplifies remarkably: 
\beq{eq:sone}
\begin{aligned}
\eta \ddot{\theta}+\alpha \dot{\theta}+\dot{\phi}& = h_3(t),\\
\eta \ddot{\phi}+\alpha \dot{\phi}-\dot{\theta}&=0. 
\end{aligned} \tag{S'}
\eeq
Although missing most of the nonlinear phenomena taking place in \eqref{eq:szero}, this new system captures its essential structural properties, which allow us to construct an effective switching field.

\section{Control construction}\label{sec:control}
Let us now introduce the following notation
\beq{eq:x}
(\theta, \dot{\theta}, \phi, \dot{\phi})=:(x_1, x_2, x_3, x_4) \mx{.}
\eeq
In this way, system \eqref{eq:sone} can be clearly written as 
\[
\dot{\ue{x}}=\uue{A} \ue{x}+\ue{f}(t),
\]
where  
\[
\uue{A}=
\begin{pmatrix}
0&1&0&0\\
0&-\alpha/\eta&0&-1/\eta\\
0&0&0&1\\
0&1/\eta&0&-\alpha/\eta
\end{pmatrix},
\qquad
\ue{f}(t)=
\begin{pmatrix}
0\\
h_3(t)/\eta\\
0\\
0
\end{pmatrix}.
\]
The formula of variation of constants gives
\beq{eq:variationconstants}
\ue{x}(t)=e^{t\uue{A}}\ue{x}(0)+\int_0^t e^{(t-s)\uue{A}}\ue{f}(s) ds,
\eeq
where 
\[
e^{t\uue{A}}
=
\begin{pmatrix}
1&I_c(t)&0&-I_s(t)\\
0&e^{-\alpha t/\eta}\cos\dfrac{t}{\eta}&0&-e^{-\alpha t/\eta}\sin\dfrac{t}{\eta}\\
0&I_s(t)&1&I_c(t)\\
0&e^{-\alpha t/\eta}\sin\dfrac{t}{\eta}&0&e^{-\alpha t/\eta}\cos\dfrac{t}{\eta}
\end{pmatrix},
\]
and
\beq{eq:icis}
\begin{aligned}
I_c(t)
& =
\frac{\eta}{1+\alpha^2}
\left[
\alpha
+
e^{-\alpha t/\eta}
\left(
\sin\frac{t}{\eta} -\alpha\cos\frac{t}{\eta}
\right)
\right]\\
I_s(t)
& =
\frac{\eta}{1+\alpha^2}
\left[
1
-
e^{-\alpha t/\eta}
\left(
\cos\frac{t}{\eta}+\alpha\sin\frac{t}{\eta}
\right)
\right].
\end{aligned}
\eeq
We now choose a class of exponentially decaying controls of the form
\beq{eq:genericcontrol}
h_3(t)=A e^{-\gamma \alpha t/\eta} \mx{,}
\eeq
where $A,\gamma >0$. By using $\ue{x}(0)=0$ and (\ref{eq:genericcontrol}) in (\ref{eq:variationconstants}), we get in particular
\beq{eq:xonexthree}
\begin{aligned}
x_1(t)
=\frac{A\eta}{1+\alpha^2}\Bigg\{
&\frac{1-e^{-\gamma\alpha t/\eta}}{\gamma}+\frac{1}{1+(\gamma-1)^2\alpha^2}
\Bigg[
\left(1+(\gamma-1)\alpha^2\right)e^{-\gamma\alpha t/\eta}\\
&\qquad-e^{-\alpha t/\eta}
\left(
(2-\gamma)\alpha\sin\frac{t}{\eta}
+\left(1+(\gamma-1)\alpha^2\right)\cos\frac{t}{\eta}
\right)
\Bigg]
\Bigg\},\\
x_3(t)
=\frac{A\eta}{1+\alpha^2}\Bigg\{
&\frac{1-e^{-\gamma\alpha t/\eta}}{\gamma\alpha}-\frac{1}{1+(\gamma-1)^2\alpha^2}
\Bigg[
(2-\gamma)\alpha e^{-\gamma\alpha t/\eta}\\
&\qquad+e^{-\alpha t/\eta}
\left(
\left(1+(\gamma-1)\alpha^2\right)\sin\frac{t}{\eta}
-(2-\gamma)\alpha\cos\frac{t}{\eta}
\right)
\Bigg]
\Bigg\}.
\end{aligned}
\eeq
Observe that, by choosing
\beq{eq:a}
A=\frac{\pi\gamma\alpha(1+\alpha^2)}{\eta},
\eeq
we obtain
\[
\lim_{t\to\infty}x_3(t)=\pi,
\qquad
\lim_{t\to\infty}x_1(t)=\pi\alpha.
\]
On the other hand, the second and fourth components read as
\[
\begin{aligned}
x_2(t)
&=\frac{\pi\gamma\alpha(1+\alpha^2)}
{\eta\left[1+(\gamma-1)^2\alpha^2\right]}
\Bigg[
e^{-\alpha t/\eta}
\left(
(\gamma-1)\alpha\cos\frac{t}{\eta}
+\sin\frac{t}{\eta}
\right)
-(\gamma-1)\alpha e^{-\gamma\alpha t/\eta}
\Bigg],\\
x_4(t)
&=\frac{\pi\gamma\alpha(1+\alpha^2)}
{\eta\left[1+(\gamma-1)^2\alpha^2\right]}
\Bigg[
e^{-\alpha t/\eta}
\left(
(\gamma-1)\alpha\sin\frac{t}{\eta}
-\cos\frac{t}{\eta}
\right)
+e^{-\gamma\alpha t/\eta}
\Bigg],
\end{aligned}
\]
from which
\[
\lim_{t\to\infty}x_2(t)
=\lim_{t\to\infty}x_4(t)=0.
\]
Hence, for small $\alpha$, the limiting state
\[
\lim_{t\to\infty}\ue{x}(t)
=(\pi\alpha,0,\pi,0)^{\mathsf T}
\]
lies close to the desired switched configuration
$(0,0,\pi,0)^{\mathsf T}$.
\begin{rem} It is easy to show that by choosing suitable nonzero initial conditions $\ue{x}(0)$, it is possible to obtain $\lim_{t \to \infty} x_1(t)=0$ exactly. However, in order to achieve this, an initial ''kick'', for instance by applying a suitable $h_2(t)$ for a short time, is necessary. The limiting state $(\pi\alpha,0,\pi,0)$ differs from the target equilibrium by an $O(\alpha)$ latitude offset.
{As shown by an energy argument in
Section~\ref{sec:relaxation}, this limiting state belongs
to the field-free basin of attraction of the reversed
equilibrium under the small-damping condition stated there.
This motivates retaining the simpler control without an
additional initial pulse.}
\end{rem}
The solution is particularly simple for $\gamma=1$, and it reads as
\beq{eq:deltaone}
\begin{aligned}
x_1(t)
&=\pi\alpha\left[
1-e^{-\alpha t/\eta}
\left(\cos\frac{t}{\eta}+\alpha\sin\frac{t}{\eta}\right)
\right],\\
x_2(t)
&=\frac{\pi\alpha(1+\alpha^2)}{\eta}
e^{-\alpha t/\eta}\sin\frac{t}{\eta},\\
x_3(t)
&=\pi\left[
1-e^{-\alpha t/\eta}
\left(
1+\alpha\sin\frac{t}{\eta}
+\alpha^2\left(1-\cos\frac{t}{\eta}\right)
\right)
\right],\\
x_4(t)
&=\frac{\pi\alpha(1+\alpha^2)}{\eta}
e^{-\alpha t/\eta}\left(1-\cos\frac{t}{\eta}\right).
\end{aligned}
\eeq

\subsection{The role of the decay parameter $\gamma$}
\label{subsec:decay-parameter}

The fact that $\gamma=1$ gives rise to solutions as simple as those obtained in \eqref{eq:deltaone}, supports the evidence that $\alpha/\eta$ is a natural rate for a control decay and it is the reason for which it has been chosen of the form \eqref{eq:genericcontrol} in the first instance. \\
 
{Increasing $\gamma$ shortens the applied-field decay time
$\eta/(\gamma\alpha)$, but does not necessarily shorten the
time required for certified confinement or relaxation.
The auxiliary solution retains an oscillatory transient
decaying on the fixed time scale $\eta/\alpha$, whose
amplitude depends on $\gamma$.}
However, it is easy to expect that this parameter cannot be increased arbitrarily. The first reason is that even the auxiliary system is less ''controlled'' for large $\gamma$ as we are about to show in this section. The second, more profound aspect, clearly related to the first one, is that the non-linear effects of the complete system could amplify those behaviours, potentially resulting in a failed switching phenomenon. This last aspect will be investigated in Section~\ref{sec:validation}.\\   
Let us start noticing that the very ''equatorial character'' of the approximated solution crucially depends on $\gamma$. In fact, from the first of (\ref{eq:xonexthree}) and (\ref{eq:a}), one gets 
\[
\begin{aligned}
x_1(t)-\pi\alpha
&=-\pi\alpha e^{-\gamma\alpha t/\eta}
+\frac{\pi\gamma\alpha}{1+(\gamma-1)^2\alpha^2}
\Bigg[
\left(1+(\gamma-1)\alpha^2\right)e^{-\gamma\alpha t/\eta}\\
&-e^{-\alpha t/\eta}
\left(
(2-\gamma)\alpha\sin\frac{t}{\eta}
+\left(1+(\gamma-1)\alpha^2\right)\cos\frac{t}{\eta}
\right)
\Bigg]\\
&=\frac{\pi\alpha(\gamma-1)(1+\alpha^2)}
{1+(\gamma-1)^2\alpha^2}e^{-\gamma\alpha t/\eta}\\
&-\frac{\pi\gamma\alpha}{1+(\gamma-1)^2\alpha^2}
e^{-\alpha t/\eta}
\left[
(2-\gamma)\alpha\sin\frac{t}{\eta}
+\left(1+(\gamma-1)\alpha^2\right)\cos\frac{t}{\eta}
\right].
\end{aligned}
\]
By using $\left|
(2-\gamma)\alpha\sin\beta
+\left(1+(\gamma-1)\alpha^2\right)\cos\beta
\right|
\leq
\sqrt{(1+\alpha^2)\left(1+(\gamma-1)^2\alpha^2\right)}$, after having taken the absolute value, we obtain
\[
|x_1(t)-\pi\alpha|
\leq
\frac{\pi\alpha|\gamma-1|(1+\alpha^2)}
{1+(\gamma-1)^2\alpha^2}
e^{-\gamma\alpha t/\eta}
+
\frac{\pi\gamma\alpha\sqrt{1+\alpha^2}}
{\sqrt{1+(\gamma-1)^2\alpha^2}}
e^{-\alpha t/\eta}.
\]

{For $\gamma\geq1$, denote the two coefficients in the
preceding bound by
\[
 C_1(\gamma)=
 \frac{\pi\alpha(\gamma-1)(1+\alpha^2)}
      {1+(\gamma-1)^2\alpha^2},
 \qquad
 C_2(\gamma)=
 \frac{\pi\gamma\alpha\sqrt{1+\alpha^2}}
      {\sqrt{1+(\gamma-1)^2\alpha^2}}.
\]
The first coefficient is non-monotone: it increases up to
$\gamma=1+\alpha^{-1}$ and decreases thereafter. The second
increases up to $\gamma=1+\alpha^{-2}$, where its maximum
is $\pi(1+\alpha^2)$.\\
However, the location of this maximum does not determine the
range in which the equatorial approximation is accurate.
When $\gamma\alpha\ll1$, the coefficients satisfy
\[
 C_1(\gamma)\simeq\pi\alpha(\gamma-1),
 \qquad C_2(\gamma)\simeq\pi\gamma\alpha,
\]
so the bound controls the latitude by a quantity of order
$\gamma\alpha$. 
For weak damping, the scale $\gamma=O(\alpha^{-1})$ at which
small latitude excursions can no longer be inferred from
the bound is much smaller than the scale
$\gamma=O(\alpha^{-2})$ at which $C_2$ attains its maximum.
For damping coefficients of order $10^{-2}$, these scales
are of order $10^2$ and $10^4$, respectively. Thus, the
large value of the maximizing decay parameter does not
justify extending the equatorial approximation to that
range; the relevant smallness requirement remains
$\gamma\alpha\ll1$.}

Thus, increasing $\gamma$ shortens the applied-field decay
time but can produce larger departures from the equator.
A quantitative restriction based on a prescribed latitude
tolerance is developed in Section~\ref{sec:relaxation}.
Moreover, since the amplitude in \eqref{eq:a} grows linearly
with $\gamma$, shorter pulses require higher peak fields.

It is also worth stressing that 
\begin{equation}\label{eq:pulse-area}
\int_0^{\infty} h_3(t;\gamma) dt =\pi (1+\alpha^2) \mx{,} 
\end{equation}
i.e. independent of $\gamma$, so that the control looks like a smoothed Dirac delta centered in $t=0$ in which $\gamma$ plays the role of smoothing parameter. In other terms, one can write, in the distributional sense
\[
h_3(\cdot;\gamma)
\xrightarrow[\gamma\to\infty]{}
\pi(1+\alpha^2)\delta_0 \,,
\]
where $\delta_0$ denotes the Dirac measure concentrated at the origin.

\section{Magnetization relaxation to equilibrium}
\label{sec:relaxation}

We define the energy without the applied-field
contribution following \eqref{eq:free energy}
\begin{equation}
 E_0(\mathbf m,\dot{\mathbf m})
 :=\frac{\eta}{2}|\dot{\mathbf m}|^2
   +\frac12\sum_{j=1}^3D_jm_j^2.
 \label{eq:relax-energy}
\end{equation}
For the field-free iLLG equation, the energy balance \eqref{eq:iLLG energy balance} yields
\[
 \frac{dE_0}{dt}=-\alpha|\dot{\mathbf m}|^2\leq0.
\]
Since $|\mathbf m|=1$, every state with $m_1=0$ satisfies
\[
 E_0\geq\frac12(D_2m_2^2+D_3m_3^2)\geq\frac{D_2}{2}.
\]
Consequently, a state with $m_1<0$ and $E_0<D_2/2$ cannot
subsequently cross the plane $m_1=0$. The energy dissipation
identity and the invariance principle then imply convergence to
$(\mathbf m,\dot{\mathbf m})=(-\mathbf e_1,\mathbf0)$, the only
field-free equilibrium in the negative hemisphere below this
energy level.

For the prescribed amplitude
$A=\gamma\pi\alpha(1+\alpha^2)/\eta$, the auxiliary solution tends to
$(\theta,\dot\theta,\phi,\dot\phi)=(\pi\alpha,0,\pi,0)$.
The corresponding resting state is
\[
 \overline{\mathbf m}
 =(-\cos(\pi\alpha),0,\sin(\pi\alpha)),
 \qquad \dot{\overline{\mathbf m}}=\mathbf0 \,.
\]
Evaluating the physical energy of the full iLLG system
at this reconstructed auxiliary state gives
\[
 E_0(\overline{\mathbf m},\mathbf0)
 =\frac{D_1}{2}
  +\frac{D_{3,1}}{2}\sin^2(\pi\alpha).
\]
Thus, if
\begin{equation}
 0<\alpha<\frac12,
 \qquad D_{3,1}\sin^2(\pi\alpha)<D_{2,1},
 \label{eq:relax-limit-condition}
\end{equation}
this state belongs to the basin of attraction of the reversed
equilibrium for the field-free full system. This statement concerns the field-free full system initialized
at $(\overline{\mathbf m},\mathbf0)$; it does not assert that
the driven full trajectory reaches this auxiliary limiting state. For fixed anisotropy
parameters, these conditions hold for all sufficiently small
$\alpha$.

The exponentially decaying field need not be removed to obtain
an energy confinement criterion. Indeed, for
$h_3(t)=A e^{-\gamma\alpha t/\eta}$,
\[
 \frac{dE_0}{dt}
 =-\alpha|\dot{\mathbf m}|^2+h_3(t)\dot m_3
 \leq\frac{h_3(t)^2}{4\alpha}.
\]
Introduce the remaining-input bound and the corrected energy
\begin{equation}
 R(t):=\frac{\eta h_3(t)^2}{8\gamma\alpha^2},
 \qquad B(t):=E_0(t)+R(t).
 \label{eq:relax-corrected-energy}
\end{equation}
Since $\dot R=-h_3^2/(4\alpha)$, the full dynamics satisfies
\begin{equation}
 \dot B(t)=-\alpha\left|
 \dot{\mathbf m}(t)-\frac{h_3(t)}{2\alpha}\mathbf e_3
 \right|^2\leq0.
 \label{eq:relax-monotonicity}
\end{equation}
Therefore, the conditions
\begin{equation}
 m_1(T)<0,
 \qquad B(T)<\frac{D_2}{2}
 \label{eq:relax-certificate}
\end{equation}
exclude backhopping for every $t\geq T$, even if the exponential
tail remains applied. They also guarantee safe field removal
at any later time. 
For a trajectory generated by the continuously applied
exponential field, define
\begin{equation}
 T_{\mathrm{sw}}
 :=
 \inf\left\{
 T\geq0:
 m_1(T)<0,\quad B(T)<\frac{D_2}{2}
 \right\},
 \label{eq:certified-switch-off-time}
\end{equation}
with $T_{\mathrm{sw}}=+\infty$ if the set is empty.
Since $B$ is nonincreasing, once the strict inequalities
hold they remain satisfied at every later time. Thus, when
$T_{\mathrm{sw}}$ is finite, field removal at any
$T>T_{\mathrm{sw}}$ is certified safe. The infimum itself
need not satisfy the strict inequalities.

Throughout the paper, $T_{\mathrm{sw}}$ denotes this
certification threshold, rather than the first crossing of
$m_1=0$, the completion of relaxation, or the earliest field
removal that could produce successful reversal.
These conditions concern the full trajectory and do not
follow from the auxiliary limiting state alone.

The auxiliary solution nevertheless provides a simple criterion
for selecting decay rates consistent with the equatorial
approximation. Choose a small angular tolerance $\delta$ such that
\begin{equation}
 0<\delta<\frac{\pi}{2},\qquad \delta\ll1,
 \qquad D_{3,1}\sin^2\delta<D_{2,1}.
 \label{eq:relax-latitude-tolerance}
\end{equation}
The last inequality ensures that a resting state with azimuth
$\pi$ and latitude at most $\delta$ lies below the barrier.
It does not replace the kinetic-energy condition during motion.
The smallness of $\delta$ controls the geometric approximation;
for example, $|\theta|\leq\delta$ implies
$|\cos^2\theta-1|\leq\sin^2\delta$.

For $\gamma\geq1$, the explicit auxiliary solution consists of
the limiting latitude $\pi\alpha$, a decaying forced term, and
a damped oscillatory term. Bounding their amplitudes gives
\begin{align}
 \sup_{t\geq0}|\theta_{\mathrm a}(t)|
 &\leq \pi\alpha
 +\pi\alpha(\gamma-1)(1+\alpha^2)
 +\pi\gamma\alpha\sqrt{1+\alpha^2}
 \notag\\
 &=\pi\alpha\left[
 \gamma\bigl(1+\alpha^2+\sqrt{1+\alpha^2}\bigr)-\alpha^2
 \right].
 \label{eq:relax-latitude-bound}
\end{align}
Consequently, a sufficient practical restriction is
\begin{equation}
 1\leq\gamma\leq
 \gamma_{\max}^{\mathrm a}(\delta,\alpha)
 :=\frac{\delta/(\pi\alpha)+\alpha^2}
 {1+\alpha^2+\sqrt{1+\alpha^2}}
 \simeq\frac{\delta}{2\pi\alpha},
 \label{eq:relax-gamma-criterion}
\end{equation}
provided the upper bound is at least $1$. This guarantees
$|\theta_{\mathrm a}(t)|\leq\delta$ for all times and quantifies
the requirement $\gamma\alpha\ll1$. Alternatively, a tolerance
$\varepsilon$ on $|\cos^2\theta_{\mathrm a}-1|$ can be imposed by
taking $\delta=\arcsin\sqrt{\varepsilon}$, subject to
\eqref{eq:relax-latitude-tolerance}. A less conservative estimate
can be obtained by maximizing the exact auxiliary latitude
numerically instead of using \eqref{eq:relax-latitude-bound}.

The bound \eqref{eq:relax-gamma-criterion} is an admissibility
criterion for the auxiliary approximation, not a theorem that
the full system switches for every rate satisfying it. In
particular, the omitted anisotropy torques can determine capture
near the saddle separating the two wells. Under
\eqref{eq:relax-limit-condition}, the auxiliary trajectory
eventually satisfies the energy and hemisphere inequalities
for every fixed $\gamma>0$, whereas the full trajectory need not.
The validation section will compare their predicted relaxation
times and examine the sharp increase of $T_{\mathrm{sw}}$ near
a finite critical decay rate, beyond which the full trajectory
can return to the original well.

An estimate of the decay rate minimizing the certified
switch-off time can also guide the control design. Introduce
\[
 s=\frac{\alpha T}{\eta},\qquad
 C=\frac{\pi^2(1+\alpha^2)^2}{8\eta},\qquad
 \Delta=\frac{D_{2,1}-D_{3,1}\sin^2(\pi\alpha)}{2}>0.
\]
For $\gamma>1$ and $\gamma\alpha\ll1$, replace the auxiliary
potential energy by its limiting value and retain the remaining
field input and the dominant oscillatory kinetic energy.
The oscillatory kinetic contribution is obtained from the
explicit auxiliary velocities in Section~\ref{sec:control}, using
$|\dot{\mathbf m}_{\mathrm a}|^2
\simeq\dot\theta_{\mathrm a}^2+\dot\phi_{\mathrm a}^2$
in the small-latitude approximation.
Neglecting the faster forced velocity and its interference with
the oscillatory component gives the approximate threshold equation
\[
 C\left[\gamma e^{-2\gamma s}
       +4\alpha^2\gamma^2e^{-2s}\right]=\Delta.
\]
The two terms describe competing effects: increasing $\gamma$
accelerates the decay of the field tail but increases the
amplitude of the inertial transient. Implicit differentiation
shows that an interior optimum satisfies
\begin{equation}
 \begin{aligned}
 C\left[\gamma e^{-2\gamma s}
       +4\alpha^2\gamma^2e^{-2s}\right]&=\Delta,\\
 (1-2\gamma s)e^{-2\gamma s}
       +8\alpha^2\gamma e^{-2s}&=0.
 \end{aligned}
 \label{eq:relax-optimal-rate}
\end{equation}
Solving these two scalar equations gives estimates of
$\gamma_{\mathrm{opt}}$ and
$T_{\min}=(\eta/\alpha)s_{\mathrm{opt}}$.
A simpler initial estimate follows by balancing the separate
decay times
\[
 T_{\mathrm{tail}}\simeq\frac{\eta}{2\alpha\gamma}
     \log\frac{C\gamma}{\Delta},\qquad
 T_{\mathrm{kin}}\simeq\frac{\eta}{2\alpha}
     \log\frac{4C\alpha^2\gamma^2}{\Delta}.
\]
Writing $L=\log(C/\Delta)$, and neglecting $\log\gamma$
relative to $L$ and $\log(4C\alpha^2/\Delta)$ relative to
$2\log\gamma$, their crossover reduces to
$2\gamma\log\gamma\simeq L$. Thus, when these logarithmic
approximations are appropriate and $L>0$,
\begin{equation}
 \gamma_{\mathrm{opt}}\simeq
 \frac{L}{2W_0(L/2)},\qquad
 L=\log\frac{\pi^2(1+\alpha^2)^2}{8\eta\Delta},
 \label{eq:relax-optimal-rate-explicit}
\end{equation}
where $W_0$ is the principal Lambert $W$ function. This explicit
estimate provides an initial guess for
\eqref{eq:relax-optimal-rate}. The chosen rate must also respect
the latitude restriction \eqref{eq:relax-gamma-criterion}; if the
estimated optimum lies outside that range, the control should
instead be optimized over the admissible interval.

Eliminating the kinetic-energy term from
\eqref{eq:relax-optimal-rate} gives
\[
 \frac{C\gamma_{\mathrm{opt}}}{2}
 \bigl(1+2\gamma_{\mathrm{opt}}s_{\mathrm{opt}}\bigr)
 e^{-2\gamma_{\mathrm{opt}}s_{\mathrm{opt}}}
 =\Delta.
\]
Consequently, the minimum certified switch-off time is estimated by
\begin{equation}
 T_{\mathrm{sw}}^{\min}\simeq
 \frac{\eta}{2\alpha\gamma_{\mathrm{opt}}}
 \left[
 -1-W_{-1}\!\left(
 -\frac{2\Delta}{eC\gamma_{\mathrm{opt}}}
 \right)
 \right],
 \label{eq:relax-minimum-time}
\end{equation}
where $W_{-1}$ denotes the real lower branch of the Lambert
$W$ function (real-valued if $0<2\Delta/(C\gamma_\mathrm{opt})\leq 1$). This estimate uses the same approximations as
\eqref{eq:relax-optimal-rate}.

An explicit estimate of the switch-off threshold can be obtained
without solving an implicit equation. Let
\[
 B_{\mathrm a}(T):=
 E_0(\mathbf m_{\mathrm a}(T),\dot{\mathbf m}_{\mathrm a}(T))
 +R(T)
\]
denote the corrected physical energy evaluated on the
reconstructed auxiliary trajectory. Replacing its potential
energy by the limiting value and neglecting the faster forced
velocity and its interference with the oscillatory component
gives the approximation
\[
 B_{\mathrm a}(T)-\frac{D_2}{2}
 \simeq-\Delta+C\gamma e^{-2\gamma\alpha T/\eta}
 +\frac{4C\alpha^2\gamma^2}{d_\gamma}
 e^{-2\alpha T/\eta},
 \qquad d_\gamma=1+(\gamma-1)^2\alpha^2.
\]
Requiring each positive contribution to be smaller than
$\Delta/2$ gives the practical design estimate
\begin{equation}
 T_{\mathrm{sw}}^{\mathrm est}(\gamma)=
 \frac{\eta}{2\alpha}
 \max\left\{
 0,\;
 \frac1\gamma\log\frac{2C\gamma}{\Delta},\;
 \log\frac{8C\alpha^2\gamma^2}{d_\gamma\Delta}
 \right\}.
 \label{eq:relax-explicit-time}
\end{equation}
The allocation of half the energy margin to each retained
contribution controls their sum within the approximate
two-exponential model. It does not bound the neglected
contributions and therefore does not provide a rigorous
upper bound on the exact auxiliary or full-system threshold.
Although the preceding argument assumes $\gamma>1$, we also
evaluate the formula at $\gamma=1$ in the numerical validation.
In that case, the forced and oscillatory velocity contributions
have the same exponential decay rate, so neglecting the former
cannot be justified by faster decay; the resulting time is
used as a heuristic design estimate.

This formula estimates the time required for both the remaining
field input and the inertial transient to become sufficiently
small. It is intended for decay rates consistent with the
auxiliary approximation and does not itself certify confinement
of the full nonlinear trajectory. A guaranteed safe removal
time is obtained by checking
\begin{equation}\label{eq:fieldfree-confinement}
 m_1(T)<0,\qquad E_0(T)<D_2/2
\end{equation}
on the full trajectory. If, more strongly,
$B(T)<D_2/2$, field removal is safe at every later time as well.
Thus, \eqref{eq:relax-explicit-time} supplies a candidate time,
while the full-system energy inequality supplies the guarantee.

For practical reference, Table~\ref{tab:control-design-formulas}
summarizes the pulse parameters, physical-unit conversions,
confinement conditions and analytical design estimates.
The table distinguishes exact relations and sufficient
full-system conditions from approximate formulas used
to guide control selection.

\begingroup

\fontsize{8}{10}\selectfont
\renewcommand{\arraystretch}{1.2}
\begin{longtable}{@{}>{\raggedright\arraybackslash}p{0.24\textwidth}
                    >{\raggedright\arraybackslash}p{0.71\textwidth}@{}}
\caption{Reference formulas for inertia-based macrospin control.
Equation references identify the corresponding formulas or
derivations in the text. Exact auxiliary relations are
distinguished from sufficient full-system conditions and
approximate design estimates. The trajectory starts from
$(\mathbf e_1,\mathbf0)$.}
\label{tab:control-design-formulas}\\
\toprule
Quantity / status / reference & Formula and interpretation\\
\midrule
\endfirsthead
\multicolumn{2}{l}{\tablename\ \thetable\ (continued)}\\
\toprule
Quantity / status / reference & Formula and interpretation\\
\midrule
\endhead
\midrule
\multicolumn{2}{r}{Continued on next page}\\
\endfoot
\bottomrule
\endlastfoot

Applied field
\par\textit{Exact prescription}
\par Eqs.~\eqref{eq:genericcontrol}, \eqref{eq:a}
& $\displaystyle
 h_3(t)=A_\gamma e^{-t/\tau_h},\qquad
 A_\gamma=\frac{\pi\gamma\alpha(1+\alpha^2)}{\eta},\qquad
 \tau_h=\frac{\eta}{\gamma\alpha}.$
\par The field is directed along $\mathbf e_3$; $\gamma>0$.
\\[5pt]

Physical units
\par\textit{Exact conversion}
\par Sec.~\ref{sec:setup}
& $\displaystyle
 \mu_0H_{3,\max}=(\mu_0M_s)A_\gamma,\qquad
 t_{\mathrm{phys}}=t_0t,\qquad
 t_0=\frac{1}{\gamma_e M_s}.$
\par Thus $\tau_{h,\mathrm{phys}}=t_0\tau_h$.
\\[5pt]

Inertial time scales
\par\textit{Auxiliary dynamics}
\par Eqs.~\eqref{eq:icis}, \eqref{eq:genericcontrol}
& $\displaystyle
 \tau_{\mathrm d}=\frac{\eta}{\alpha},\qquad
 \omega_{\mathrm n}=\frac1\eta,\qquad
 \frac{\tau_h}{\tau_{\mathrm d}}=\frac1\gamma.$
\par $\tau_{\mathrm d}$ is the velocity-amplitude damping time;
$\omega_{\mathrm n}$ is the angular frequency of the auxiliary mode.
\\[5pt]

Pulse area
\par\textit{Exact}
\par Eq.~\eqref{eq:pulse-area}
& $\displaystyle
 \int_0^\infty h_3(t)\,dt=\pi(1+\alpha^2).$
\par The normalized pulse area is independent of $\gamma$ and $\eta$.
\\[5pt]

Auxiliary limiting state
\par\textit{Exact within the auxiliary model}
\par Eqs.~\eqref{eq:xonexthree}, \eqref{eq:a}
& $\displaystyle
 (\theta_{\mathrm a},\dot\theta_{\mathrm a},
  \phi_{\mathrm a},\dot\phi_{\mathrm a})
 \longrightarrow(\pi\alpha,0,\pi,0).$
\par This is an auxiliary prediction, not the full-system equilibrium.
\\[5pt]

Energy and barrier
\par\textit{Full macrospin model}
\par Eq.~\eqref{eq:relax-energy}
& $\displaystyle
 E_0=\frac\eta2|\dot{\mathbf m}|^2
       +\frac12\sum_{j=1}^3D_jm_j^2,\qquad
 E_{\mathrm b}=\frac{D_2}{2}.$
\par Every state with $m_1=0$ has $E_0\geq E_{\mathrm b}$.
\\[5pt]

Safe removal at $T$
\par\textit{Sufficient condition}
\par Eq.~\eqref{eq:fieldfree-confinement}
& $\displaystyle m_1(T)<0,\qquad E_0(T)<E_{\mathrm b}.$
\par Evaluate on the full trajectory. Removal at this time is
followed by field-free confinement and relaxation to $-\mathbf e_1$.
\\[5pt]

Confinement with a field tail
\par\textit{Sufficient condition}
\par Eqs.~\eqref{eq:relax-corrected-energy},
\eqref{eq:relax-certificate}
& $\displaystyle
 \begin{aligned}
 R(t)&=\frac{\eta h_3(t)^2}{8\gamma\alpha^2}
      =C\gamma e^{-2\gamma\alpha t/\eta},\\
 B(t)&=E_0(t)+R(t),\qquad
 m_1(T)<0,\quad B(T)<E_{\mathrm b}.
 \end{aligned}$
\par Excludes subsequent backhopping with the tail retained;
also certifies removal at any later time.
\\[5pt]

Latitude-based rate limit
\par\textit{Sufficient auxiliary bound}
\par Eq.~\eqref{eq:relax-gamma-criterion}
& $\displaystyle
 1\leq\gamma\leq\gamma_{\max}^{\mathrm a}
 =\frac{\delta/(\pi\alpha)+\alpha^2}
 {1+\alpha^2+\sqrt{1+\alpha^2}}
 \simeq\frac{\delta}{2\pi\alpha}.$
\par Guarantees $\sup_t|\theta_{\mathrm a}(t)|\leq\delta$,
provided the upper bound is at least $1$.
Choose $0<\delta\ll1$ with
$D_{3,1}\sin^2\delta<D_{2,1}$.
This bounds latitude, not full-system capture.
\\[5pt]

Candidate confinement time
\par\textit{Approximate}
\par Eq.~\eqref{eq:relax-explicit-time}
& $\displaystyle
 T_{\mathrm{sw}}^{\mathrm{est}}(\gamma)
 =\frac\eta{2\alpha}\max\left\{
 0,\frac1\gamma\log\frac{2C\gamma}{\Delta},
 \log\frac{8C\alpha^2\gamma^2}{d_\gamma\Delta}
 \right\}.$
\par An explicit design estimate; verification of a full-system
confinement condition is still required.
\\[5pt]

Optimal rate and time
\par\textit{Coupled approximation}
\par Eq.~\eqref{eq:relax-optimal-rate}
& $\displaystyle
 \begin{aligned}
 C[\gamma e^{-2\gamma s}+4\alpha^2\gamma^2e^{-2s}]&=\Delta,\\
 (1-2\gamma s)e^{-2\gamma s}+8\alpha^2\gamma e^{-2s}&=0.
 \end{aligned}$
\par Solve for an interior candidate
$(\gamma_{\mathrm{opt}},s_{\mathrm{opt}})$,
with $T_{\min}=(\eta/\alpha)s_{\mathrm{opt}}$.
The approximation assumes $\gamma>1$ and $\gamma\alpha\ll1$.
\\[5pt]

Explicit rate estimate
\par\textit{Rough initial guess}
\par Eq.~\eqref{eq:relax-optimal-rate-explicit}
& $\displaystyle
 \gamma_{\mathrm{LW}}\simeq\frac{L}{2W_0(L/2)}.$
\par Requires $L>0$ and the logarithmic approximations used
in the derivation. It does not minimize
$T_{\mathrm{sw}}^{\mathrm{est}}(\gamma)$ exactly.
\\[5pt]

Minimum-time estimate
\par\textit{Approximate}
\par Eq.~\eqref{eq:relax-minimum-time}
& $\displaystyle
 T_{\min}\simeq
 \frac{\eta}{2\alpha\gamma_{\mathrm{opt}}}
 \left[-1-W_{-1}\!\left(
 -\frac{2\Delta}{eC\gamma_{\mathrm{opt}}}
 \right)\right].$
\par Evaluated at the coupled candidate optimal rate.
The real branch requires
$0<2\Delta/(C\gamma_{\mathrm{opt}})\leq1$.
Substituting $\gamma_{\mathrm{LW}}$ gives a further rough estimate.
\\
\end{longtable}
Here $W_0$ and $W_{-1}$ denote the principal and lower real
branches of the Lambert $W$ function, respectively. The
design estimates apply to the auxiliary approximation and
must respect any imposed latitude tolerance. Only the stated
full-system energy inequalities provide confinement guarantees;
their direct use for spatially nonuniform micromagnetic states
is not justified by the macrospin argument.

\endgroup

\section{Validation and numerical experiments}
\label{sec:validation}

We first assess the proposed control in the macrospin setting,
comparing the auxiliary solution with the full nonlinear iLLG
dynamics. The comparison distinguishes the initial driven motion
from the subsequent relaxation and tests the energy criterion
developed in Section~\ref{sec:relaxation}.

\subsection{Macrospin validation}
\label{subsec:macrospin-validation}

Throughout the macrospin experiments we use the dimensionless
parameters
\begin{equation}
 \alpha=0.023,\qquad \eta=0.0105,\qquad
 (D_1,D_2,D_3)=(-0.3,0,1),
 \label{eq:validation-parameters}
\end{equation}
so that $D_{2,1}=0.3$, $D_{3,1}=1.3$, and the field-free
energy barrier is $D_2/2=0$. The magnetization is initially
at rest at the original easy-axis minimum:
\[
 \mathbf m(0)=\mathbf e_1,\qquad
 \dot{\mathbf m}(0)=\mathbf0,
\]
or, equivalently,
$(\theta,\phi,\dot\theta,\dot\phi)(0)=(0,0,0,0)$.
The same applied field is used in the full and auxiliary systems:
\begin{equation}
 \mathbf h(t)=\left(0,0,
 \frac{\pi\gamma\alpha(1+\alpha^2)}{\eta}
 e^{-\gamma\alpha t/\eta}\right),\qquad t\geq0.
 \label{eq:validation-field}
\end{equation}
The field remains applied throughout each simulation according
to this decay law; it is not switched off at $t=4$.

The spherical-coordinate equations are integrated using
\texttt{ode45}. For the comparisons on $0\leq t\leq4$,
the relative and absolute tolerances are $10^{-9}$ and $10^{-11}$,
respectively, and the maximum time step is
\[
 \Delta t_{\max}
 =\min\left\{10^{-3},\frac{\eta}{20},
                    \frac{\eta}{20\gamma\alpha}\right\}.
\]
For the long-time experiment, the full nonlinear system is
integrated to $t=10^3$, with relative and absolute tolerances
$10^{-8}$ and $10^{-10}$. On $0\leq t\leq5$, its maximum
time step is additionally restricted to
$\min\{\eta/5,\eta/(20\gamma\alpha)\}$ to resolve the initial
transient. All times and energies reported below are dimensionless.

\begin{figure*}[t!]\begin{center}
  		\begin{minipage}[c][1\width]{0.49\textwidth}
			\hspace{0pt}
			{\begin{overpic}[width=1\textwidth]{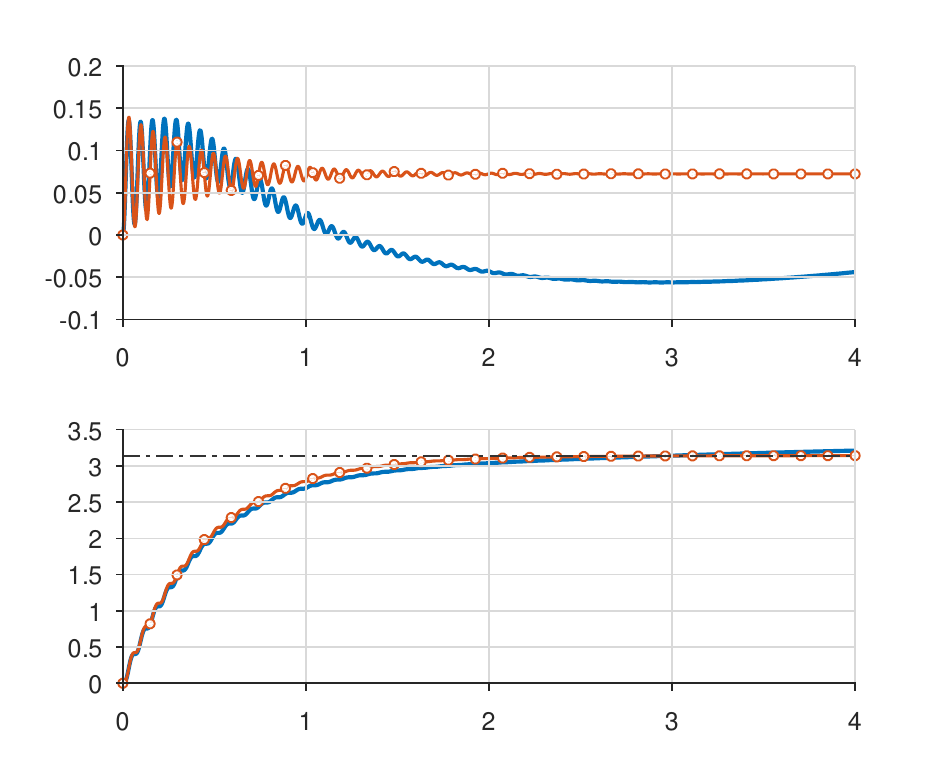}
					\put(-3,78){\footnotesize (a)}
           			\put(-3,58){\footnotesize $\theta(t)$}  
                   	\put(-3,23){\footnotesize $\phi(t)$}  
                   	\put(90,41){\footnotesize $t$}  
                   	\put(90,2){\footnotesize $t$}  
			\end{overpic}}
		\end{minipage}	
		\begin{minipage}[c][1\width]{0.49\textwidth}
			\hspace{5pt} 
			{\begin{overpic}[width=\textwidth]{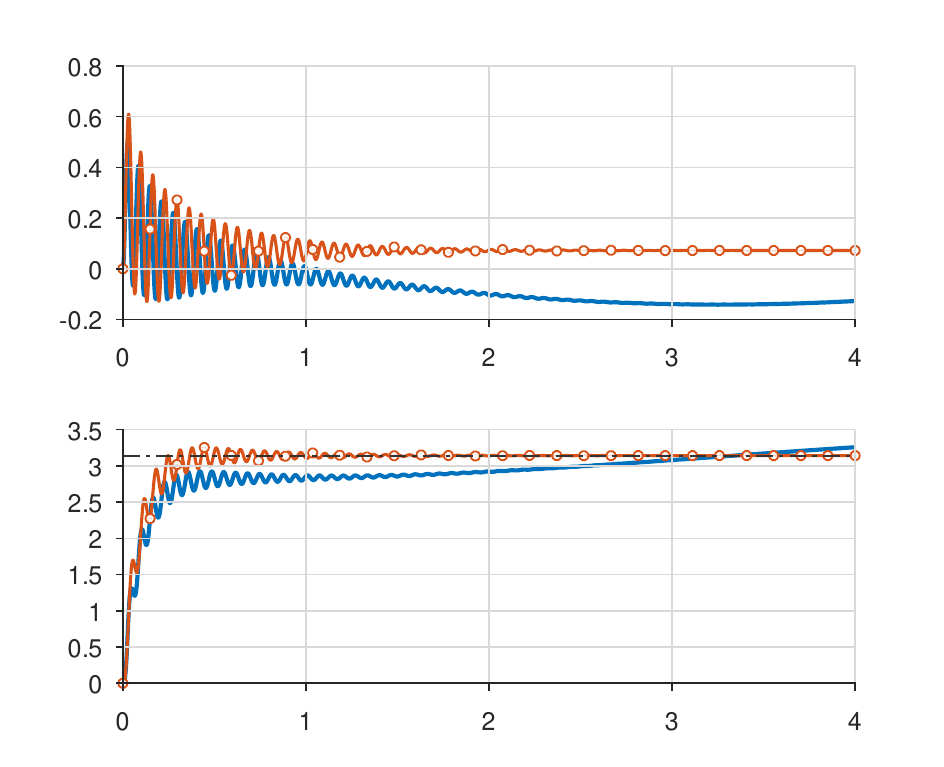}
					\put(-3,78){\footnotesize (b)}
           			\put(-3,58){\footnotesize $\theta(t)$}  
                   	\put(-3,23){\footnotesize $\phi(t)$}  
                   	\put(90,41){\footnotesize $t$}  
                   	\put(90,2){\footnotesize $t$}  
			\end{overpic}}
		\end{minipage}\\
		\vspace{-30pt}
  	\begin{minipage}[c][1\width]{0.49\textwidth}
			\hspace{5pt} 
			{\begin{overpic}[width=\textwidth]{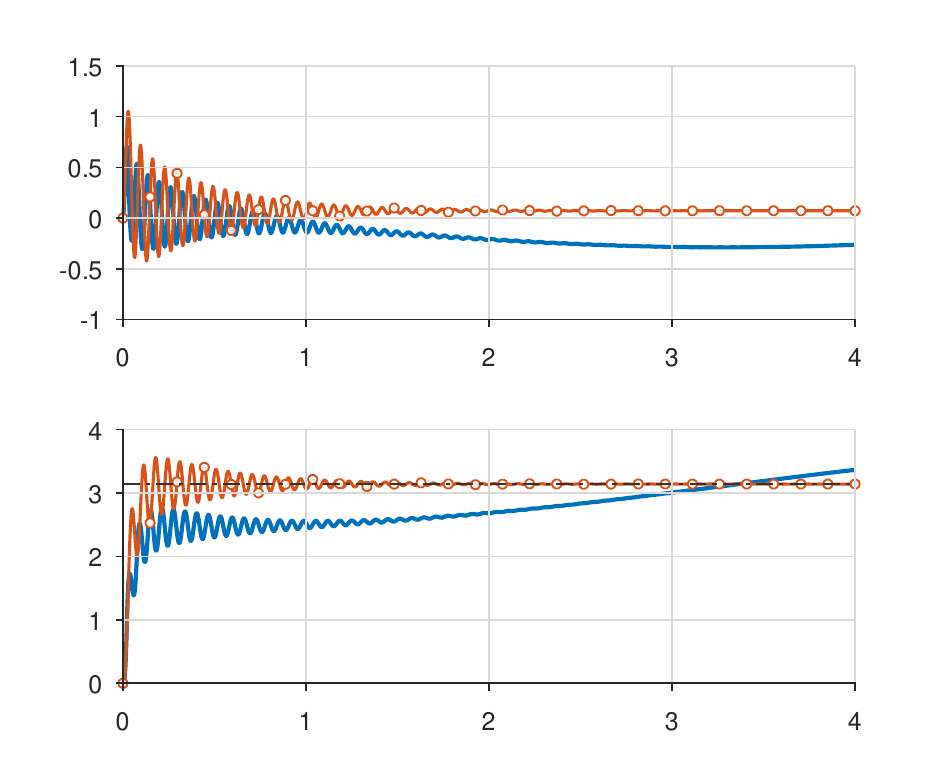}
					\put(-3,78){\footnotesize (c)}
           			\put(-3,58){\footnotesize $\theta(t)$}  
                   	\put(-3,23){\footnotesize $\phi(t)$}  
                   	\put(90,41){\footnotesize $t$}  
                   	\put(90,2){\footnotesize $t$}               
			\end{overpic}}
		\end{minipage}	
    \begin{minipage}[c][1\width]{0.50\textwidth}
			\hspace{5pt}
			{\begin{overpic}[width=\textwidth]{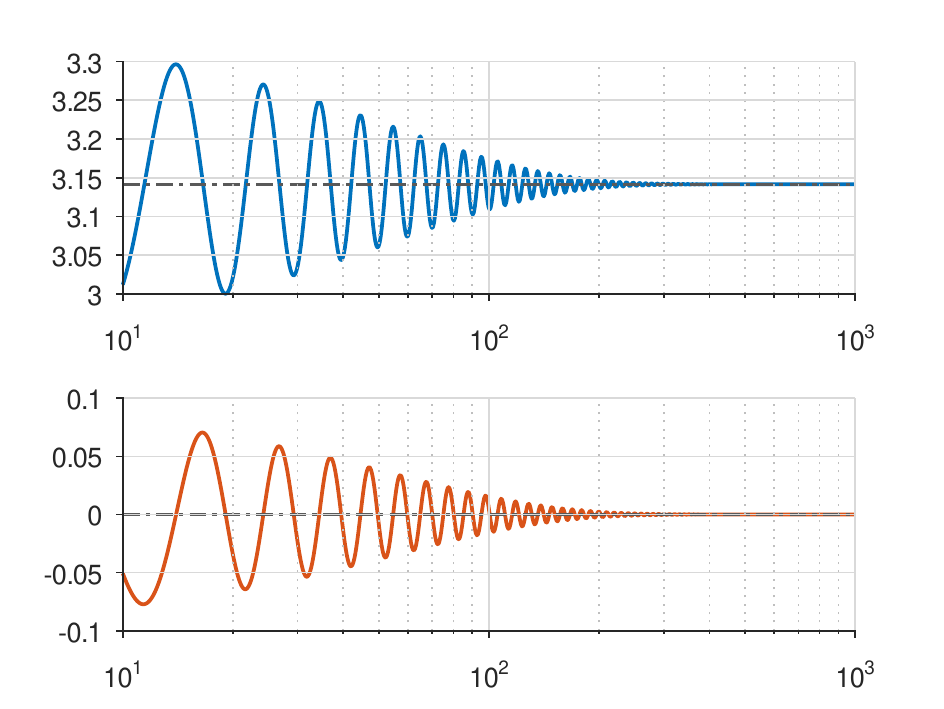}

                    \put(-3,78){\footnotesize (d)}
           			\put(-3,58){\footnotesize $\phi(t)$}  
                   	\put(-3,23){\footnotesize $\theta(t)$}  
                   	\put(70,37.5){\footnotesize $t$}  
                   	\put(70,2){\footnotesize $t$}             
 			\end{overpic}}
		\end{minipage}\\
	\end{center}
	\vspace{-30pt}
	\caption{In panels $(a)$, $(b)$ and $(c)$, comparison between the controlled solutions (blue) alongside the solution of the auxiliary system (orange - circles), for three different values of $\gamma$. More precisely, we have set $\gamma=1$ in $(a)$, $\gamma=5$ in $(b)$ and $\gamma=10$ in $(c)$. The dashed-dotted line denotes the target asymptote $\phi=\pi$.\\
    Panel $(d)$ shows the long-time evolution of the full nonlinear
system under the exponentially decaying control with $\gamma=10$,
illustrating relaxation towards the reversed equilibrium
$(\theta,\phi)=(0,\pi)$. }
	\label{fig:validation1} 
	\vspace{0pt}
\end{figure*}

Figure~\ref{fig:validation1}(a)--(c) compares the full solution
(blue) with the auxiliary solution (orange with circles) for
$\gamma=1,5,10$, respectively. For $\gamma=1$, the auxiliary
model closely reproduces the initial azimuthal rotation and
the rapid latitude oscillations. Increasing $\gamma$ produces
a more concentrated excitation and a faster initial advance
of the azimuth, but also larger latitude excursions and more
pronounced differences between the two models. In particular,
the maximum auxiliary latitudes are approximately $0.1395$,
$0.6104$, and $1.0552$ radians for $\gamma=1,5,10$.
Thus, the latter cases illustrate the limitations of interpreting
the auxiliary motion as uniformly quasi-equatorial, even though
the resulting controls remain effective in these full-system
experiments.

The distinction between the two models is also apparent after
the fast transient. The auxiliary solution approaches
$(\theta_{\mathrm a},\phi_{\mathrm a})=(\pi\alpha,\pi)$,
with $\pi\alpha\simeq0.07226$, whereas the full system retains
the anisotropy torques and continues to relax. Consequently,
the residual latitude and azimuthal deviations of the full
trajectory at $t=4$ should not be interpreted as its limiting
state. Figure~\ref{fig:validation1}(d), for $\gamma=10$,
extends the observation interval to $t=10^3$ and shows the
subsequent approach to $(\theta,\phi)=(0,\pi)$, corresponding
to the reversed equilibrium $\mathbf m=-\mathbf e_1$.
The applied field is negligible during this late evolution,
but no separate field-removal event is imposed.

The energy criterion gives a quantitative complement to this
geometric picture. For the present parameters,
\[
 E_0(t)=\frac{\eta}{2}
 \left(\dot\theta^2+\cos^2\theta\,\dot\phi^2\right)
 -0.15\cos^2\theta\cos^2\phi+0.5\sin^2\theta.
\]
Evaluating $B=E_0+R$ from
\eqref{eq:relax-corrected-energy} along the full trajectories
gives the values in Table~\ref{tab:validation-confinement}.
All three cases satisfy $m_1(4)<0$ and $B(4)<0$.
By \eqref{eq:relax-monotonicity}, subsequent backhopping is
therefore excluded, including the effect of the remaining
exponential tail. The time $T=4$ is an observation time at
which confinement is certified, not the earliest safe time
and not a prescribed switch-off time. 

Table~\ref{tab:validation-confinement} compares the explicit
estimate \eqref{eq:relax-explicit-time} with the thresholds
computed from the exact auxiliary solution and the full
nonlinear dynamics. It also reports the full-system confinement
checks at $T=4$. The exponential field remains applied throughout
these computations; the reported thresholds identify the onset
of the corresponding confinement conditions, rather than an
imposed field-removal event.

\begin{table}[tbp]
 \centering
 \begin{tabular}{cccccc}
 \hline
 $\gamma$
 & $T_{\mathrm{sw}}^{\mathrm est}$
 & $T_{\mathrm{sw}}^{\mathrm a}$
 & $T_{\mathrm{sw}}$
 & $m_1(4)$
 & $B(4)$ \\
 \hline
 $1$  & $1.684695$ & $1.529881$ & $1.532696$
      & $-0.996629$ & $-0.148009$ \\
 $5$  & $1.011832$ & $0.849321$ & $0.814867$
      & $-0.985330$ & $-0.137457$ \\
 $10$ & $1.320615$ & $1.165742$ & $1.135488$
      & $-0.941375$ & $-0.098587$ \\
 \hline
 \end{tabular}
 \caption{Switch-off threshold estimates and energy confinement
 checks for the controls in Figure~\ref{fig:validation1}.
 $T_{\mathrm{sw}}^{\mathrm est}$ is the explicit analytical estimate,
 $T_{\mathrm{sw}}^{\mathrm a}$ is obtained from the exact
 auxiliary solution, and $T_{\mathrm{sw}}$ is the full-system
 certified threshold. The last two columns refer to the full
 trajectory at $T=4$, with barrier $D_2/2=0$.
 All times are dimensionless.}
 \label{tab:validation-confinement}
\end{table}

The explicit estimate exceeds the full-system threshold by
approximately $10\%$, $24\%$, and $16\%$ for $\gamma=1,5,10$,
respectively. The field-tail contribution determines
$T_{\mathrm{sw}}^{\mathrm est}$ for $\gamma=1$, whereas the inertial
relaxation contribution determines it for $\gamma=5,10$.
This change reflects the competition underlying the
nonmonotone dependence of the certified time on $\gamma$.
The exact auxiliary thresholds are closer to the full-system
values, with relative differences below $4.3\%$ in these tests.

For all three cases, $T_{\mathrm{sw}}^{\mathrm est}>T_{\mathrm{sw}}$.
Consequently, the full-system energy certificate guarantees
safe field removal at the estimated times and at every later
time. Moreover, $m_1(4)<0$ and $B(4)<0$ independently confirm
confinement at the observation time $T=4$, including the
remaining exponential tail. The explicit formula therefore
provides a conservative design estimate in these examples;
this numerical observation does not establish its
conservativeness for arbitrary decay rates or parameters.

These results support the use of the auxiliary solution for
control design while showing why successful reversal must be
assessed using the full dynamics. A shorter applied pulse
does not by itself imply faster relaxation or a smaller
certified switch-off time.

\begin{figure}[tbp]
 \centering
 \includegraphics[width=\textwidth]{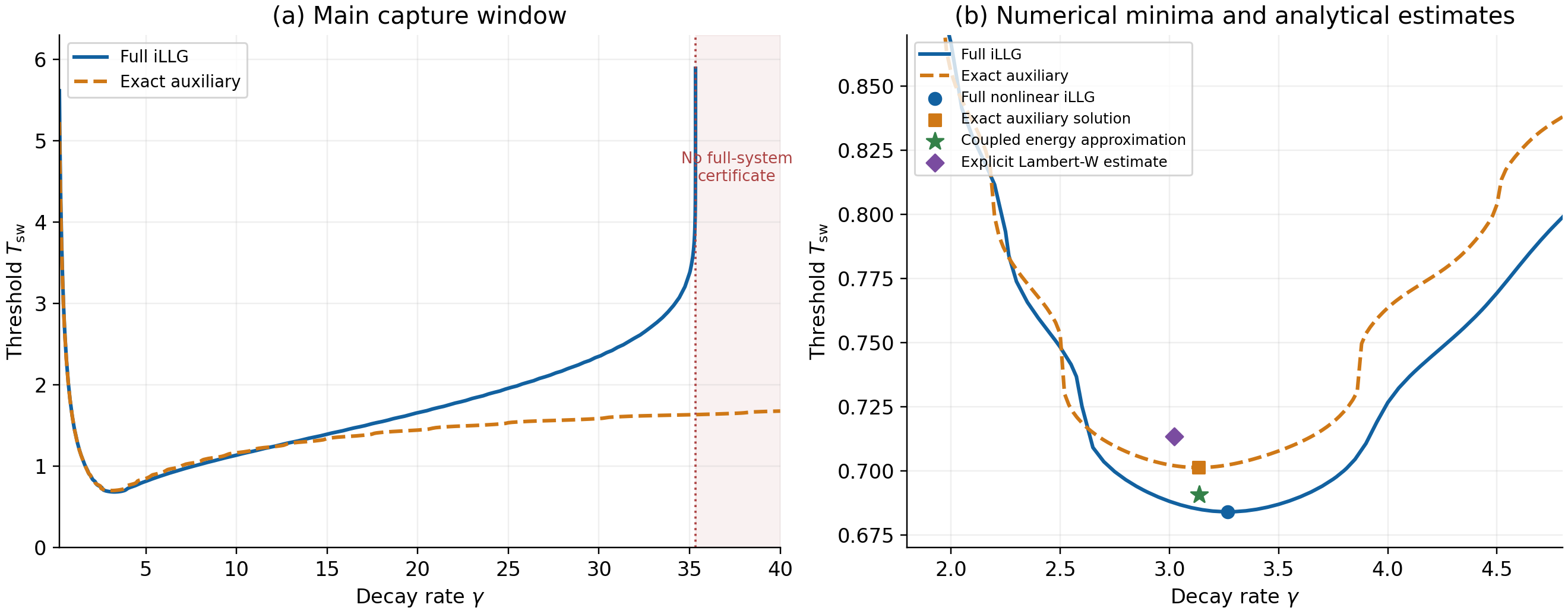}
 \caption{Certified switch-off thresholds and control-design
 estimates for the same exponential pulse family.
 (a) Full nonlinear iLLG (blue) and exact auxiliary solution
 (orange dashed) on the main capture window. The vertical
 dotted line marks the estimated upper boundary
 $\gamma_c\simeq35.326377$; the shaded region indicates sampled
 rates with capture in the original well.
 (b) Enlarged view of the minima. Blue circles and orange squares
 identify the numerical full and auxiliary minima; the green
 star and purple diamond give the coupled energy approximation
 and explicit Lambert-$W$ estimate, respectively.
 Parameter values are given in \eqref{eq:validation-parameters}.}
 \label{fig:switching-time-design}
\end{figure}

\paragraph{Decay-rate selection and minimum certified time.}
Figure~\ref{fig:switching-time-design} compares the onset of
certified confinement for the full and auxiliary trajectories
as the decay rate varies. For the full dynamics,
$T_{\mathrm{sw}}$ is obtained from the downward crossing
$B(T)=0$, provided $m_1(T)<0$. For the auxiliary trajectory,
the same physical energy is evaluated on the reconstructed
magnetization, and $T_{\mathrm{sw}}^{\mathrm a}$ is the last
boundary after which both $B_{\mathrm a}<0$ and
$m_{\mathrm a,1}<0$ remain satisfied. This distinction is
needed because $B_{\mathrm a}$ need not be monotone.
The reported thresholds denote the onset; strict confinement
inequalities hold at slightly later times.

The two predictions agree closely for the moderate rates
illustrated in Figure~\ref{fig:validation1}, but their
discrepancy increases near the upper boundary of the main
capture window. The full threshold grows sharply as
$\gamma$ approaches $\gamma_c\simeq35.326377$, while the
auxiliary threshold remains finite and varies smoothly.
The sampled trajectories just beyond this boundary become
trapped in the original hemisphere and admit no persistent
certificate of reversed confinement. This observation does
not exclude an earlier, deliberately timed field removal.
The sharp increase is consistent with a prolonged passage
near the anisotropy saddle, whose capture dynamics is absent
from the auxiliary model. The numerical scan alone does not
constitute a proof of divergence at $\gamma_c$.

The minima found on the main switching branch are summarized
in Table~\ref{tab:optimal-decay}. The coupled energy estimate
is obtained from \eqref{eq:relax-optimal-rate}, whereas the
explicit estimate uses \eqref{eq:relax-optimal-rate-explicit}
together with \eqref{eq:relax-minimum-time}.
The explicit pair provides a rough initial estimate; the
coupled equations retain more of the energy balance.
In particular, the coupled approximation predicts the optimal
decay rate within about $4.1\%$ and the minimum time within
about $1.0\%$ of the full-system values. The exact auxiliary
solution also locates the optimal rate closely, despite its
different late-time relaxation. These are unconstrained
optima: imposing a specified latitude tolerance requires
optimization over the corresponding admissible interval.

\begin{table}[tbp]
 \centering
 \begin{tabular}{lcc}
 \hline
 Method & $\gamma_{\mathrm{opt}}$ & $T_{\min}$\\
 \hline
 Full nonlinear iLLG          & $3.26781$ & $0.683965$\\
 Exact auxiliary solution    & $3.13278$ & $0.701226$\\
 Coupled energy approximation & $3.13508$ & $0.690740$\\
 Explicit Lambert-$W$ estimate & $3.02277$ & $0.713359$\\
 \hline
 \end{tabular}
 \caption{Numerical minima of the certified switch-off threshold
 and analytical design estimates for the parameters in
 \eqref{eq:validation-parameters}.}
 \label{tab:optimal-decay}
\end{table}

The backhopping examples and refinement of the full-system
minimum use the DOP853 integrator with relative and absolute
tolerances $2\times10^{-13}$ and $2\times10^{-15}$.
Maximum time steps are $0.002$ for the backhopping examples
and $0.005$ for refinement of the minimum. The auxiliary
thresholds are evaluated from the explicit solution, with
crossings resolved numerically and a decaying envelope used
to bound the subsequent motion.

\begin{figure}[tbp]
 \centering
 \includegraphics[width=\textwidth]{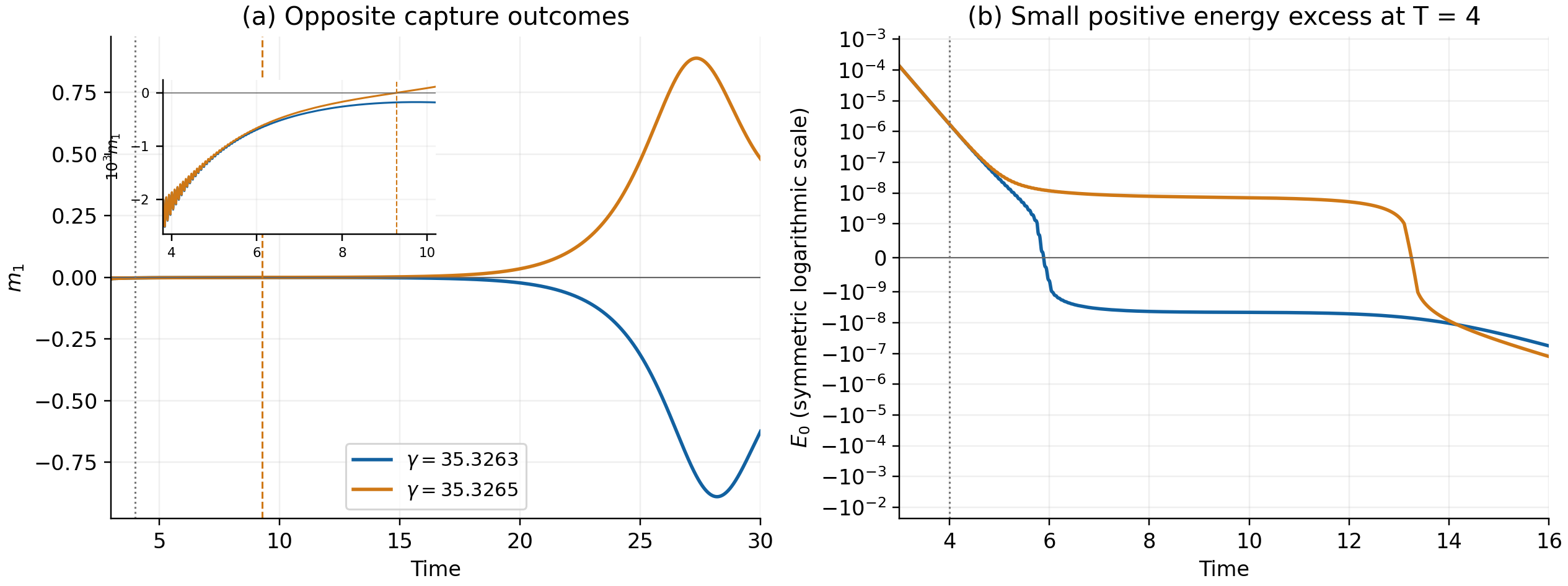}
 \caption{Full nonlinear trajectories for two nearby exponential
 pulses with $\gamma=35.3263$ (blue) and $35.3265$ (orange).
 (a) Evolution of $m_1$; the inset resolves the small
 magnetization components near the backhopping event.
 The orange dashed line marks $t\simeq9.28420$.
 (b) Energy relative to the barrier $D_2/2=0$, using a symmetric
 logarithmic scale with a linear region $|E_0|\leq10^{-9}$.
 The grey dotted lines mark $T=4$. Both trajectories have
 a small positive energy at this time, but only the orange
 trajectory subsequently returns to the original hemisphere.}
 \label{fig:backhopping-pulse}
\end{figure}

\paragraph{Backhopping near the energy barrier.}
The decay-rate comparison identifies a critical region in which
small changes in $\gamma$ can alter the eventual capture outcome.
To illustrate this sensitivity and test the significance of the
strict energy inequality, we consider two nearby rates,
$\gamma=35.3263$ and $\gamma=35.3265$, on opposite sides of the
estimated boundary $\gamma_c$. The initial conditions and physical
parameters are unchanged. Figure~\ref{fig:backhopping-pulse}
shows that the two trajectories eventually become confined
to opposite hemispheres. At $T=4$,
both trajectories are in the reversed hemisphere, but slightly
above the barrier:
\[
\begin{array}{c|cc}
 \gamma & m_1(4) & E_0(4)\\ \hline
 35.3263 & -0.00191611 & 1.72561\times10^{-6}\\
 35.3265 & -0.00190539 & 1.73853\times10^{-6}.
\end{array}
\]
The remaining-input bound satisfies $R(4)<6\times10^{-266}$,
so the failure of $B(4)<0$ reflects the actual energy excess,
rather than an overestimate of the field tail.

For $\gamma=35.3265$, backhopping occurs at
$t\simeq9.28420$, when $E_0\simeq7.29\times10^{-9}>0$.
The energy subsequently becomes negative in the original
hemisphere, excluding another crossing. By contrast, for
$\gamma=35.3263$, the corrected energy becomes negative at
$t\simeq5.88727$ while $m_1<0$, certifying confinement to the
reversed hemisphere. The field is not removed in either
experiment. Thus, a small positive violation of the energy
condition can permit backhopping, but does not necessarily
cause it: the condition is sufficient, not necessary, for
successful reversal.

\subsection{Inertial micromagnetic simulations}
\label{subsec:micromagnetic-validation}

To assess the applicability of the proposed control beyond the
macrospin approximation, we consider spatially resolved inertial
micromagnetic simulations. Following the setting of
Ref.~\cite{FdS26}, the magnetic sample is a permalloy disk with
diameter $1\,\mu\mathrm m$ and thickness
$5\,\mathrm{nm}$. Its small aspect ratio
approximates the thin-film geometry underlying the macrospin
model, while allowing edge effects and spatially nonuniform
magnetization dynamics. The effective field includes exchange,
magnetostatic, uniaxial anisotropy and applied-field contributions.

We use $\mu_0M_s=1\,\mathrm T$ and exchange stiffness
$A_{\mathrm{ex}}=13\,\mathrm{pJ\,m^{-1}}$, together with the
dimensionless damping and inertial parameters
$\alpha=0.023$ and $\eta=0.0105$. The uniaxial easy axis is
parallel to $\mathbf e_1$. Writing its energy density as
$K_1(1-m_1^2)$, the choice $\frac{2K_1}{\mu_0M_s^2}=D_{2,1}=0.3$
gives $K_1\simeq1.19\times10^5\,\mathrm{J\,m^{-3}}$.
For an ideal thin film, this reproduces the macrospin anisotropy
differences $D_{2,1}=0.3$ and $D_{3,1}=1.3$.
The finite disk introduces corrections through its nonlocal
magnetostatic field, which is computed in the micromagnetic
model rather than replaced by the ideal demagnetizing factors.

The disk is discretized using cubic cells of side
$5\,\mathrm{nm}$, corresponding to a $200\times200\times1$
bounding grid. The cell size is smaller than the exchange
length $\ell_{\mathrm{ex}}= \sqrt{\frac{2A_{\mathrm{ex}}}{\mu_0M_s^2}} \simeq5.71\,\mathrm{nm}$.
%
Time integration uses the midpoint geometric method for
inertial magnetization dynamics described in Ref.~\cite{DAQUINO2024112874} and implemented in the code MaGICo~\cite{MaGICo},
with physical time step $25\,\mathrm{fs}$.
The initial configuration is obtained by relaxing the sample
in zero applied field from a uniformly magnetized state along
$\mathbf e_1$. During this preparatory stage, a large damping
$\alpha=1$ is used until the normalized
effective torque is below $10^{-5}$. The switching stage then
starts from the relaxed configuration with zero inertial
velocity and with the physical damping restored.

A spatially uniform field is applied along $\mathbf e_3$.
In the same dimensionless normalization as the macrospin model,
its prescribed decay is given by eq.\eqref{eq:validation-field}.
We compare $\gamma=1,5,10$ and the optimal rate selected from
the macrospin study, $\gamma_{\mathrm{opt}}\simeq3.02$ computed from eq.\eqref{eq:relax-optimal-rate-explicit}.
For the chosen parameters, the peak applied field is
$\mu_0H_{3,\max}=\mu_0M_s\pi\gamma\alpha(1+\alpha^2)/\eta
\simeq6.885\gamma\,\mathrm T$, corresponding to
$6.89$, $34.43$ and $68.85\,\mathrm T$ for $\gamma=1,5,10$,
respectively, and $20.81\,\mathrm T$ for the analytical
choice $\gamma_{\mathrm{opt}}\simeq3.023$ from eq.\eqref{eq:relax-optimal-rate-explicit}.
The cases with $\gamma>1$ therefore require particularly
high peak fields, which pose a substantial experimental
challenge. Their micromagnetic simulations are included
primarily for comparison with the macrospin analysis and
to assess the predicted dependence on the decay rate,
rather than as immediately accessible experimental protocols.

\begin{figure}[t!]
 \centering
 \includegraphics[width=1\linewidth]{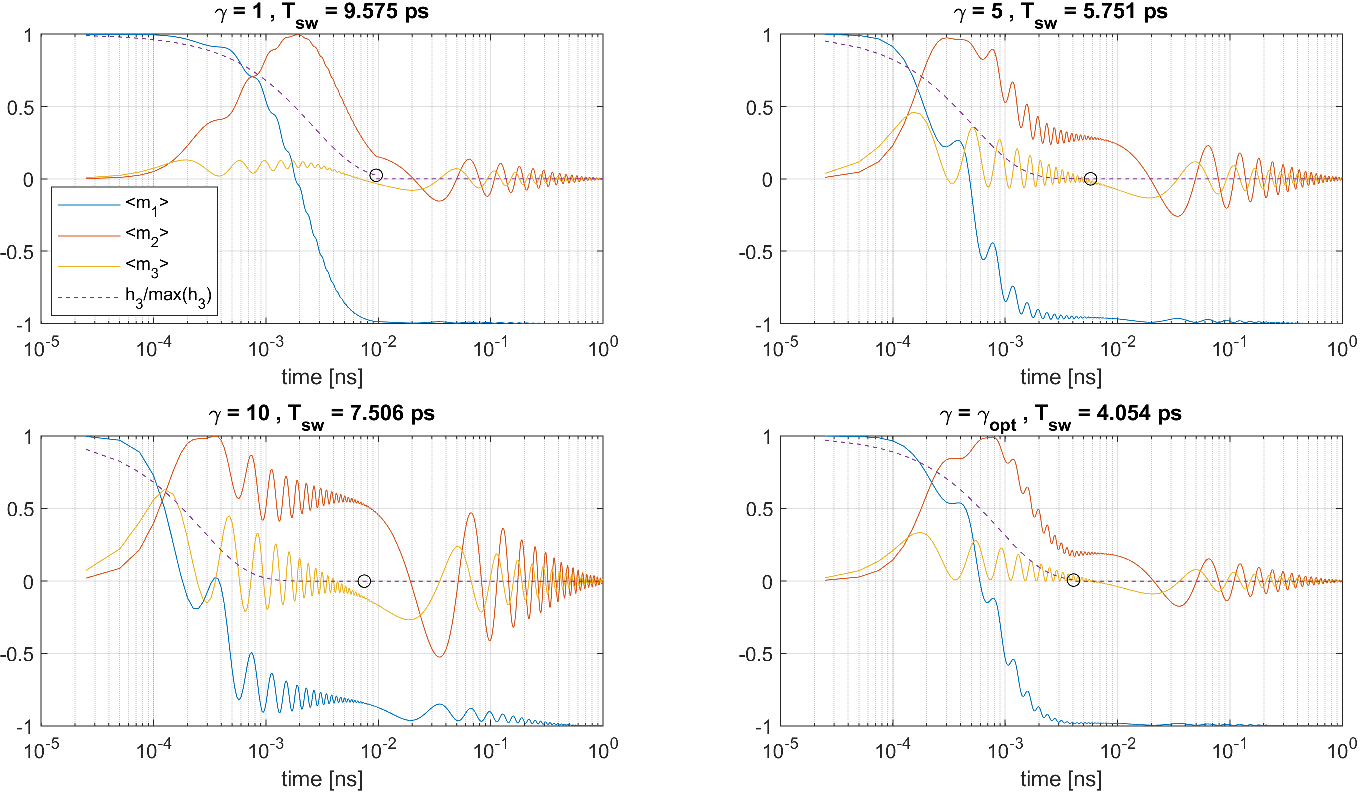}
  \caption{Inertial micromagnetic switching and relaxation under
the decaying control for $\gamma=1,5,10$ and
$\gamma_{\mathrm{opt}}$, simulated up to $1\,\mathrm{ns}$.
Each panel reports the spatially averaged magnetization
components and the field normalized by its own maximum,
$h_3(t)/\max_{t\geq0}h_3(t)$, against physical time on a
logarithmic axis. The circles mark the analytical reference
times annotated above the panels: $9.575$, $5.751$, $7.506$
and $4.054\,\mathrm{ps}$, respectively. All trajectories remain
in the reversed hemisphere after these times and approach the
reversed equilibrium through damped transverse oscillations.
The markers indicate confinement estimates, not completion
of relaxation.}
 \label{fig:mumag}
\end{figure}

Figure~\ref{fig:mumag} shows reversal of the spatially averaged
longitudinal magnetization for all four decay rates. Extending
the simulations to $t=1\,\mathrm{ns}$ resolves both the initial
switching transient and the subsequent damped relaxation.
In every case, $\langle m_1\rangle$ remains negative after the
marked analytical time and approaches $-1$, while the transverse
averages tend towards zero. No subsequent backhopping is
observed over the simulated interval.

The dependence on $\gamma$ reflects the competition identified
in the macrospin analysis. For $\gamma=1$, reversal proceeds
with a relatively smooth longitudinal evolution and a small
out-of-plane component. Increasing the rate to $\gamma=5$
and $10$ advances the initial rotation but excites larger
transverse oscillations. The nonmonotone longitudinal transient
and the oscillations of $\langle m_2\rangle$ and
$\langle m_3\rangle$ are particularly pronounced for
$\gamma=10$. The field is normalized by its own maximum in
each panel.

The logarithmic time axis distinguishes the rapid initial
oscillations from the slower damped motion that persists after
the applied field has become negligible. This later relaxation
is most pronounced for $\gamma=10$, which retains substantial
transverse excursions on the scale of tens to hundreds of
picoseconds. The extended integration shows that these
oscillations decay rather than lead to a return to the original
hemisphere. By the end of the simulated interval, all four
trajectories are close to the reversed equilibrium, with
strongly attenuated transverse motion.

The annotated analytical times are $9.575\,\mathrm{ps}$,
$5.751\,\mathrm{ps}$, and $7.506\,\mathrm{ps}$ computed from \eqref{eq:relax-explicit-time} for
$\gamma=1,5,10$, respectively, and $4.054\,\mathrm{ps}$ from \eqref{eq:relax-minimum-time}
for the panel labelled $\gamma=\gamma_{\mathrm{opt}}$.
The optimized case provides the smallest annotated reference
time among those shown while maintaining reversal throughout
the subsequent relaxation. In particular, the most rapidly
decaying excitation, $\gamma=10$, does not give the shortest
predicted confinement time and produces the largest late
transverse oscillations. The marked times therefore identify
design estimates for confinement, rather than completion of
relaxation or the first zero crossing of the longitudinal
magnetization.

\begin{figure}[tbp]
 \centering
 \includegraphics[width=0.6\linewidth]{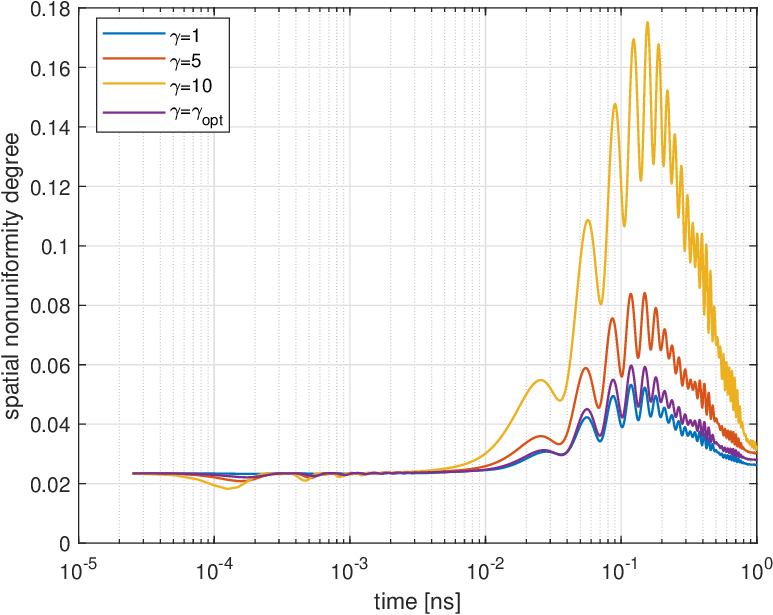}
 \caption{Spatial nonuniformity degree
 $\mathcal U(t)=\sqrt{1-|\langle\mathbf m\rangle(t)|^2}$
 for the inertial micromagnetic trajectories in
 Figure~\ref{fig:mumag}: $\gamma=1$ (blue), $5$ (orange),
 $10$ (yellow), and $\gamma_{\mathrm{opt}}$ (purple).
 Time is shown in nanoseconds on a logarithmic axis.
 Nonuniformity remains small during the initial switching
 stage and exhibits a delayed transient increase during
 relaxation, strongest for $\gamma=10$. The optimized control
 produces a peak close to that for $\gamma=1$.
 All curves decrease towards small residual values by
 $t=1\,\mathrm{ns}$.}
 \label{fig:mumag-uniformity}
\end{figure}

The spatial coherence of the reversal is quantified in
Figure~\ref{fig:mumag-uniformity} using the nonuniformity degree
\begin{equation}
 \mathcal U(t):=
 \sqrt{1-|\langle\mathbf m\rangle(t)|^2}
 =\left\langle
 |\mathbf m-\langle\mathbf m\rangle|^2
 \right\rangle^{1/2}.
 \label{eq:micromagnetic-nonuniformity}
\end{equation}
The equality follows from the pointwise constraint
$|\mathbf m|=1$. Thus, $\mathcal U$ is the root-mean-square
deviation from the spatially averaged magnetization and
vanishes for a perfectly uniform configuration. Its small
initial value, approximately $0.024$, reflects the slight
nonuniformity of the relaxed starting state.

For all four decay rates, $\mathcal U$ remains close to its
initial value during the first few picoseconds, encompassing
the initial reversal and the marked analytical confinement
times. The largest spatial variations develop later, during
the damped relaxation. The transient peaks occur on the scale
of $0.1$--$0.2\,\mathrm{ns}$ and increase markedly with the
decay rate: the approximate maxima are $0.055$, $0.085$ and
$0.175$ for $\gamma=1,5,10$, respectively. The optimized
control gives a maximum of about $0.060$, close to that for
$\gamma=1$ and substantially below those for $\gamma=5,10$.
The optimized rate therefore combines the short analytical
confinement time with comparatively weak spatial nonuniformity.

Even the largest observed peak, $\mathcal U\simeq0.175$,
corresponds to
$|\langle\mathbf m\rangle|=\sqrt{1-\mathcal U^2}\simeq0.985$.
The magnetization consequently remains strongly aligned on
average, supporting a predominantly coherent description
while resolving finite spatial deviations from the macrospin
motion. For every rate, the nonuniformity subsequently decreases
to a small value close to its initial level by $1\,\mathrm{ns}$.
Together with Figure~\ref{fig:mumag}, this shows that the
transient growth of spatial variations is followed by relaxation
towards a nearly uniform reversed state, without observed
backhopping over the simulated interval.

The global measure $\mathcal U$ does not identify the spatial
structure or location of the deviations, and does not exclude
localized modes. Moreover, the macrospin energy barrier cannot
be applied directly to a nonuniform magnetization field.
These results therefore provide numerical corroboration of
the macrospin control strategy rather than a direct extension
of its energy confinement theorem to the micromagnetic system.

\section{Conclusions}
\label{sec:conclusions}

We have constructed a family of exponentially decaying magnetic
fields for inertial magnetization reversal from an initial
easy-axis state at rest. The control is designed through an
explicitly resolvable auxiliary model, with magnetic inertia
directly determining the pulse amplitude and temporal scale.
In the parameter regimes validated numerically, the prescribed
excitation produces reversal followed by relaxation without
requiring field interruption at a particular phase of the motion.

A central advantage of the proposed strategy is its ability
to overcome the demanding switch-off timing associated with
conventional precessional switching. The energy analysis
provides a sufficient confinement criterion that accounts for
both inertial kinetic energy and the remaining work of the
applied field. Once this criterion holds in the reversed
hemisphere, subsequent backhopping is excluded: the exponential
tail may remain applied, or the field may be removed at any
later time. Thus, a narrow, phase-sensitive interruption window
is replaced by a sufficient lower time threshold for safe
field removal.

The auxiliary solution provides practical guidance for selecting
the control parameters. The latitude criterion quantifies a
regime consistent with the equatorial approximation, while the
approximate energy balance identifies the competition between
field decay and inertial relaxation. This balance predicts an
intermediate decay rate minimizing the certified confinement
time, in close agreement with the full nonlinear simulations.
The near-critical examples show that the upper boundary of
successful capture must be assessed using the full nonlinear
dynamics. For moderate decay rates, however, the analytical
estimates agree closely with the full-system results,
providing useful tools for the design of the switching process.

The inertial micromagnetic simulations further support the
effectiveness of the strategy beyond the macrospin approximation.
For all cases considered, reversal is maintained throughout
the simulated interval, while the spatially averaged
magnetization approaches the reversed equilibrium. The initial
switching is predominantly coherent, and the spatial
nonuniformity that develops during subsequent relaxation
eventually decreases. The analytically selected decay rate
combines a short predicted confinement time with comparatively
weak spatial nonuniformity.

These results establish a route towards exploiting magnetic
inertia as a resource for ultrafast magnetization control,
combining rapid reversal with reduced sensitivity to field
termination timing. The high peak fields required by the
investigated pulses, particularly at larger decay rates,
remain a substantial experimental and technological challenge.
Nevertheless, the explicit control construction and its energy
confinement criterion provide a theoretical basis for further
development and experimental assessment of inertia-based
switching protocols.

\subsection*{Acknowledgements} 
The numerical simulations and the corresponding plots have been performed with GNU Octave ~\cite{oct}. 

\bibliographystyle{alpha}
\bibliography{cqlsII.bib}

\end{document}